\documentclass[12pt]{amsart}
\usepackage{amsmath,amsthm,amsfonts,amssymb}
\begin{document} 
\newcommand{\A}{{\mathbb A}}
\newcommand{\B}{{\mathbb B}}
\newcommand{\C}{{\mathbb C}}
\newcommand{\N}{{\mathbb N}}
\newcommand{\Q}{{\mathbb Q}}
\newcommand{\Z}{{\mathbb Z}}
\renewcommand{\P}{{\mathbb P}}
\newcommand{\R}{{\mathbb R}}
\newcommand{\rc}{\subset}
\newcommand{\rank}{\mathop{rank}}
\newcommand{\ord}{\mathop{ord}{}}
\newcommand{\trace}{\mathop{tr}}
\newcommand{\dimc}{\mathop{dim}_{\C}}
\newcommand{\Lie}{\mathop{Lie}}
\newcommand{\Spec}{\mathop{Spec}}
\newcommand{\Auto}{\mathop{{\rm Aut}_{\mathcal O}}}
\newcommand{\alg}[1]{{\mathbf #1}}
\newtheorem*{definition}{Definition}
\newtheorem*{claim}{Claim}
\newtheorem{corollary}{Corollary}
\newtheorem*{Conjecture}{Conjecture}
\newtheorem*{SpecAss}{Special Assumptions}
\newtheorem{example}{Example}
\newtheorem*{remark}{Remark}
\newtheorem*{observation}{Observation}
\newtheorem*{fact}{Fact}
\newtheorem*{remarks}{Remarks}
\newtheorem{lemma}{Lemma}
\newtheorem{proposition}{Proposition}
\newtheorem{theorem}{Theorem}
\title{%
On meromorphic functions which are Brody curves
}
\author {J\"org Winkelmann}
\begin{abstract}
We discuss meromorphic functions on the complex plane which are Brody
curves regarded as holomorphic maps to $\P_1$, i.e., which have
bounded
spherical derivative.
\end{abstract}
\subjclass{}%
%
\address{%
J\"org Winkelmann \\
Mathematisches Institut\\
Universit\"at Bayreuth\\
Universit\"atsstra\ss e 30\\
D-95447 Bayreuth\\
Germany\\
}
\email{jwinkel@member.ams.org\newline\indent{\itshape Webpage: }%
http://btm8x5.mat.uni-bayreuth.de/\~{ }winkelmann/
}
\maketitle
We discuss meromorphic functions on $\C$ which are Brody
curves regarded as holomorphic curves from $\C$ to $\P_1(\C)$.
In other words, those meromorphic functions for which $||f'||$
calculated with respect to the euclidean metric on $\C$ and the
Fubini-Study-metric on $\P_1$ is bounded.

In concrete terms this means:
\[
\lim\sup \frac{|f'(z)|}{1+|f(z)|^2} < +\infty.
\]

This number is also called ``spherical derivative'' (\cite{CH}).
We follow the established notation and denote this spherical
derivative as
\[
f^\#(z)=\frac{|f'(z)|}{1+|f(z)|^2}.
\]

Brody curves in general have been studied since the seminal paper
of Brody (\cite{B}) illuminated their relevance for hyperbolicity questions.

The ``spherical derivative'' $f^\#$ has been studied before Brody
proved his
theorem, especially by Clunie, Hayman, Lehto and Virtanen
(\cite{CH},\cite{L},\cite{LV}).

In particular, Clunie and Hayman showed that a meromorphic Brody function is of
order at most two and a Brody entire function is of order at most one.
They also demonstrated the existence of Brody curves with certain
divisors with high multiplicities as prescribed zero divisors.

Easy examples for Brody functions are polynomials and the exponential
function.

As noted by Tsukamoto (\cite{T}), elliptic functions are also Brody:
If $\Gamma$ is a lattice in $\C$, and $f$ is a $\Gamma$-periodic 
meromorphic function, then $f^\#$ is continuous and $\Gamma$-invariant.
Therefore the existence of a compact fundamental region for $\Gamma$
implies that $f^\#$ is bounded.

Here we will show that Brody functions remain Brody after multiplying
with certain rational functions. We investigate two special classes
of meromorphic functions, namely $R(z)e^z+Q(z)$ where $R,Q$ are 
rational functions and $e^z+e^{\lambda z}$ with $\lambda\in\C$.
For these functions we determine completely which ones are Brody.

This provides us with some surprising examples: If $f,g$ are Brody
neither $f+g$, nor $fg$ need to be Brody. The Brody condition
is not closed: The meromorphic function $f_t(z)=e^z+\frac{z}{zt+1}$
is Brody if and only if $t\ne 0$.

In particular in view of the results of Clunie and Hayman on the
order it might be natural to assume that every divisor of sufficiently
slow growth can be realized as the zero divisor of an entire Brody
function. We show that this is not the case, there are effective
divisors of arbitrary slow growth which can not be realized as
zero divisors of a Brody function.%
\footnote{A related result has been obtained much earlier by Lehto (\cite{L}).}

\section{Basic properties}
We start with the observation
\[
f^\#=\left|\frac{f'}{f}\right|\cdot h(f)
\]
where
\[
h(w)=\frac{|w|}{1+|w|^2}.
\]
Since $h(z)\le\frac{1}{2}$ for all $z\in\C$, we may deduce:
{\em If $f'/f$ is bounded, then $f$ is Brody.}
In fact one has the following slightly stronger statement:
\begin{lemma}
Let $f$ be an entire function with
\[
\limsup_{z\to\infty}\left|\frac{f'(z)}{f(z)}\right|<+\infty.
\]
Then $f$ is Brody.
\end{lemma}
\begin{proof}
There is a compact subset $K\subset\C$ and a constant $C>0$
such that $|f(z)/f'(z)|<C$ for all $z\not\in K$.
Then
\[
\sup_{z\in\C}f^\#(z)\le\max\left\{\max_{z\in K}f^\#(z),C\right\}
\]
\end{proof}
\begin{corollary}
The exponential function $z\mapsto e^z$ is Brody.
\end{corollary}
Concerning rational functions the Brody property follows from
the observation below:
\begin{lemma}
Let $R$ be a non-constant rational function on $\C$.

Then $\lim_{z\to\infty}\frac{R'(z)}{R(z)}=0$.
\end{lemma}
\begin{proof}
Write $R$ as a quotient of polynomials: $R=P/Q$.
Then
\[
\frac{R'}{R}=\frac{P'Q-PQ'}{PQ}
\]
and the assertion follows from the fact that
$\deg(P'Q-PQ')<\deg(PQ)$.
\end{proof}
\begin{corollary}
Rational functions have the Brody property.
\end{corollary}

\begin{proposition}
Assume that $f$ is a Brody entire function, that $R$ is a rational
function and that $\alpha\in\C^*$, $\beta\in\C$.

Then $z\mapsto F(z)=R(f(\alpha z+\beta))$ is Brody, too.
\end{proposition}
\begin{proof}
Each rational function $R$ defines a continuous self-map of $\P_1$.
The compactness of $\P_1$ implies that for each rational function $R$
there exists a constant $C$ such that the operator norm of 
the differential map $DR:T\P_1\to T\P_1$ is bounded by $C$.
Then $\sup F^\#\le C|\alpha|\sup f^\#$.
\end{proof}
\begin{corollary}
The trigonometric functions $\sin$, $\cos$, $\tan$, $\sinh$,
$\cosh$.. are all Brody.
\end{corollary}
\begin{proof}
These functions can be expressed in the form $R(e^z)$ or $R(e^{iz})$
for some rational function $R\in\C(X)$.
\end{proof}

\section{Products with rational functions}
\begin{lemma}
Let $f:\C\to\P_1$ be Brody and let $R$ be a rational function
with $R(\infty)\not\in\{0,\infty\}$.

Then $g=Rf$ is Brody.
\end{lemma}

\begin{proof}
We have 
\[
g^{\#}(z)=\frac{g'(z)}{g(z)}h(g(z))
\]
with $h(w)=\frac{|w|}{1+|w|^2}$.

We will need the following auxiliary fact:

{\em Claim: Let $\lambda>1$. Then $h(\lambda w)\le\lambda h(w)$
and $h(\lambda^{-1}w)\le\lambda w$ for all $w\in\C$.}

This claim is rather immediate:
\[
h(\lambda w)=\frac{|\lambda w|}{1+|\lambda w|^2}
\le\frac{|\lambda w|}{1+|w|^2}=\lambda h(|w|)
\]
and similarily for the second inequality.

Now write the rational function $R$ as quotient of two polynomials:
$R(z)=\frac{P(z)}{Q(z)}$.

Then
\[
\frac{g'}{g}=\frac{f'}{f}+\frac{P'}{P}-\frac{Q'}{Q}
\]
and (using $\lim  (P'/P)=\lim (Q'/Q)=0$) we may deduce that there
are constants $R_1,C_1>0$ such that
\[
\left|\frac{g'(z)}{g(z)}-\frac{f'(z)}{f(z)}\right|
\le C_1
\]
for complex numbers $z$ with $|z|>R_1$.

By assumption $\lim_{z\to\infty}R(z)\in\C^*$. Therefore there are
constants $\lambda>1$ and $R_2>0$ such that
\[
\lambda^{-1}<|R(z)|<\lambda
\]
whenever $|z|>R_2$.

Then
\[
h(g(z))=h(R(z)f(z))\le \lambda h(f(z))
\]
if $|z|>R_2$.

As a consequence, we have
\begin{align*}
g^{\#}(z)=\frac{|g'(z)|}{|g(z)|}h(g(z)
&\le \left( \frac {|f'(z)|}{|f(z)|} +C_1\right)\lambda h(f(z)) \\
&\le \lambda f^{\#}(z) + C_1\lambda h(f(z)) \\
&\le \lambda \left( f^{\#}(z)+\frac{C_1}{2}\right)
\end{align*}
if $|z|>\max\{R_1,R_2\}$.

Hence the assertion.
\end{proof}

It is really necessary to assume $R(\infty)\ne 0,\infty$.
For instance, consider the functions  $f_1(z)=e^z+1$ and
$f_2(z)=1/f_1(z)$.
Both are Brody functions, but proposition~\ref{case-1}
 in the next section implies that neither
$g_1(z)=zf_1(z)$ nor $g_2(z)=\frac{1}{z}f_2(z)$ is Brody.

\section{Case study: $R(z)e^z+Q(z)$}

We will now discuss certain sums.
\begin{proposition}\label{case-1}
Let $R,Q$ be rational functions.
Then $f(z)=R(z)e^z+Q(z)$ is a Brody function if and only if
one of the following conditions is fulfilled:
\begin{enumerate}
\item $R\equiv 0$ or
\item $Q(\infty)\ne\infty$.
\end{enumerate}
\end{proposition}
\begin{proof}
If $R\equiv 0$, then $f=Q$ is rational and therefore Brody.

Assume now that $Q(\infty)\ne\infty$ and $R\not\equiv 0$. 
Then there is a rational function
$Q_0$ with $Q_0(\infty)\not\in\{0,\infty\}$ and a complex number $c$
with $Q(z)=Q_0(z)+c$.
Then
\[
f(z)=Q_0(z)\left(\frac{R(z)}{Q_0(z)}e^z+1\right)+c
\]
and $f$ is Brody if and only if $g(z)=\frac{R(z)}{Q_0(z)}e^z$ is Brody.
Now $S=R/Q_0$ is rational and $S\not\equiv 0$. Therefore
\[
\lim_{z\to\infty} \frac{g'(z)}{g(z)}=
\lim_{z\to\infty} \frac{S(z)+S'(z)}{S(z)}=1.
\]
But this implies
\[
\limsup_{z\to\infty}g^\#(z)\le 1\cdot \frac{1}{2}<+\infty.
\]
Hence $g$ is Brody and therefore $f$ is Brody, too.

It remains to show that $f$ is not Brody whenever $Q(\infty)=\infty$
and $R\not\equiv 0$.
In this case $Q(z)=Q_0(z)z^n$ for some rational function $Q_0$ with
$Q_0(\infty)\not\in\{0,\infty\}$ and $n\in\N$.
Define $S=R/Q_0$. Then $S$ is rational.
Considering the equality
\[
f(z)=Q_0(z)\left(S(z)e^z+z^n\right)
\]
is Brody if and only if $g(z)=S(z)e^z+z^n$ is Brody.
Since $R\not\equiv 0$ together with $Q_0(\infty)\ne 0$ implies
$S\not\equiv 0$, the entire function $g$ is transcendental.
Therefore there is a sequence of complex numbers $z_k$ with
$\lim_{k\to\infty}|z_k|=\infty$
and
$\lim_{k\to\infty}g(z_k)=0$.
Now
\[
g'(z)=(S(z)+S'(z))e^z+nz^{n-1}=\left(S(z)+S'(z)\right)
\frac{g(z)-z^n}{S(z)}+nz^{n-1}
\]
which implies
\[
g'(z)=z^n\left(\frac{S(z)+S'(z)}{S(z)}\left(\frac{g(z)}{z^n}-1\right)
+\frac{n}{z}\right).
\]
We recall that $\lim_{z\to\infty}\frac{S+S'}{S}=1$,
$\lim_{k\to\infty}z_k=\infty$ and
$\lim_{k\to\infty}g(z_k)=0$.
Combined, these facts yield
\[
\lim_{k\to\infty}g'(z_k)=\infty
\]
Together with $\lim_{k\to\infty}g(z_k)=0$ this implies
\[
\lim_{k\to\infty}g^\#(z_k)=\infty
\]
Thus $g$ is not Brody, which in turn implies that $f$ is not Brody.
\end{proof}

\section{Case study: $e^z+e^{\lambda z}$}
\begin{proposition}\label{ex-re}
The entire function $f(z)=e^z+e^{\lambda z}$ is Brody if and only
if $\lambda\in\R$.
\end{proposition}
\begin{proof}
{\em Case 1)} The cases $\lambda\in\{0,1,-1\}$ are trivial.

{\em Case 2)} 
Let $\lambda\not\in\R$. 
Then $f(z)=0$ iff $e^{z(1-\lambda)}=-1$. Hence the set of all zeros
of $f$ is given as the set of all
\[
a_k=\frac{(2k+1)\pi i}{1-\lambda}
\]
with $k$ running through the set of integers $\Z$.
Now $f(a_k)=0$ implies $f^\#(a_k)=|f'(a_k)|$
and
\[
f'(a_k)=(\lambda-1)e^\frac{(2k+1)\pi i \lambda}{1-\lambda}
\]
which is unbounded because $\lambda\not\in\R$ implies
that 
\[
\left|e^\frac{2\pi i \lambda}{1-\lambda}\right|\ne 1.
\]
Therefore $\sup_kf^\#(a_k)=+\infty$ and consequently $f$ is not Brody.

{\em Case 3)}
Assume $\lambda>0$. 
Observe that $f(z)=g(\lambda z)$ for 
$g(w)=e^w+e^{\frac{1}{\lambda}w}$. 
Hence we may assume without loss of generality that
$0<\lambda\le1$. 
Since $f(z)=2e^z$ if $\lambda=1$, it suffices to consider the
case $0<\lambda<1$. We choose
\[
C=\frac{\log 2}{1-\lambda}.
\]
Then
\[
|e^{\lambda z}|\le \frac{1}{2}|e^z|
\]
for all $z\in\C$ with $\Re(z)\ge C$.
It follows that
\[
f^\#(z)
=
\frac{|e^z+\lambda e^{\lambda z}|}{1+|e^z+e^{\lambda z}|^2}
\le\frac{3/2|e^z|}{1+1/4|e^z|^2}\le 
\frac{3/2|e^z|}{1/4|e^z|^2}=6e^{-\Re z}\le 6e^{-C}
\]
for all $z\in\C$ with $\Re(z)\ge C$.

On the other hand
\[
f^\#(z)\le|f'(z)|\le e^C+\lambda e^{\lambda C}
\]
for all $z$ with $\Re(z)\le C$.
Thus 
\[
f^\#(z)\le\max\{6e^{-C},e^C+\lambda e^{\lambda C}\}\ \forall z\in\C
\]
and consequently $f$ is Brody.

{\em Case 4)}
We assume $\lambda<0$.
It suffices to consider the case where $-1<\lambda<0$.

As before, we define
\[
C=\frac{\log 2}{1-\lambda}.
\]
and obtain
\[
|e^{\lambda z}|\le\frac{1}{2}|e^z|
\]
and therefore
\[
f^\#(z)\le 6e^{-C}
\]
for all $z\in\C$ with $\Re(z)\ge C$.

Next we observe that
the condition $\Re(z)\le -C$ implies
\[
\left|e^z\right|\le\frac12\left| e^{\lambda z}\right|.
\]
As a consequence
\[
f^\#(z)
\le\frac{3/2|e^{\lambda z}|}{1+1/4|e^{\lambda z}|^2}\le 
\frac{3/2|e^{\lambda z}|}{1/4|e^{\lambda z}|^2}=6e^{-\Re\lambda z}\le
6e^{\lambda C}
\]
for all $z$ with $\Re(z)\le -C$.

Finally we observe that
\[
f^\#(z)\le |f'(z)|\le|e^z|+|\lambda e^{\lambda z}|\le
e^C+|\lambda|e^{-\lambda C}
\]
for all $z$ with $-C\le\Re(z)\le C$.
\end{proof}
\section{Divisors of slow growth}
\begin{theorem}\label{div-slow-gro}
Let $D$ be the divisor defined by
$D=\{a_k:k\in\N\}$ (all multiplicities being one). 
Assume that there is a number $\lambda>1$ such that
\begin{enumerate}
\item
$|a_{k+1}|>\lambda|a_k|>0$ for all $k$, and
\item
The origin $0$ is contained in the interior of the convex hull
of the set of accumulation points of
the sequence $\frac{a_k}{|a_k|}$.
\end{enumerate}

Then there does not exist any Brody function with $D$ as zero divisor.
\end{theorem}

\begin{proof}
The first condition implies the absolute convergence
of
\[
\sum_k\frac{1}{|a_k|}
\]
which implies the convergence of
\[
\sum_k\log\left(1-\frac{1}{|a_k|}\right)
\]
which in turn implies the convergence of 
\[
F(z)=\prod_k\left(1-\frac{z}{a_k}\right)
\]
Thus an arbitrary entire function $f$ with divisor $D$ can be written
as
\[
f(z)=F(z)e^{g(z)}
\]
where $F$ is defined as above while $g$ is an entire function.

\begin{claim}
There is a constant $C>0$ (depending only on $\lambda$) such that
\[
\left|\prod_{k\ge n}\left(1-\frac{p}{a_k}\right)\right|\ge C
\]
for every $n\in\N$ and $p\in\{z\in\C:|p|\sqrt{\lambda}<|a_n|\}$.
\end{claim}

To prove the claim, we observe that 
\[
\left|\prod_{k\ge n}\left(1-\frac{p}{a_k}\right)\right|\ge
C=\prod_{l\ge 0}\left(1-\lambda^{-(l+1/2)}\right).
\]

We emphasize that $C$ is independent of $n$.

We need a second claim.
\begin{claim}
For every $K>0$ and $m\in\N$ there is a natural number $N\ge m$
(depending on $K$ and $m$)
such that
\[
\left|\prod_{k=m}^n\left(1-\frac{p}{a_k}\right)\right|
\ge K
\]
for all $n\ge N$ and $p$ with $|p|\ge\sqrt{\lambda}|a_n|$.
\end{claim}

To prove the second claim it suffices to note that
\[
\left|\prod_{k=m}^n\left(1-\frac{p}{a_k}\right)\right|
\ge
\prod_{l=0}^{n-m}\left(\lambda^{1/2+l}-1\right)
\]
and
\[
\lim_{s\to\infty}\prod_{l=0}^{s}\left(\lambda^{1/2+l}-1\right)=\infty.
\]

Next we choose for each $n\in\N$ a complex number $p_n$ such that
\[
\sqrt \lambda|a_n|< |p_n|<\frac{|a_{n+1}|}{\sqrt \lambda}
\]
Due to the two above claims we know that there is a number $M\in\N$
such that
\[
|F(p_n)|=\left|\prod_{k=1}^\infty\left(1-\frac{p_n}{a_k}\right)\right|\ge
 1
\]
for all $n\ge M$.

If we set $r_n=\sqrt{|a_{n+1}a_n|}$,
we can now deduce:
\[
M_{\exp(g)}(r_n)=\max_{|z|=r_n}|g(z)|\le M_f(r_n).
\]
Hence, if $f$ is Brody and therefore $T_f(r)=O(r)$, we can
deduce that $T_{\exp(g)}(r)=O(r)$, which
in turn implies that $g$ is a polynomial of degree 1,
i.e.~affine-linear.
Hence condition two implies that there
is a subsequence $a_{n_j}$ with
\[
\left| \exp(g(a_{n_j}))\right|\ge 1 \ \forall j
\]
because $|e^w|=e^{\Re(w)}$.

Next we will calculate $|f^{\#}(a_{n_j})|$.
Since $f(a_k)=F(a_k)=0$ for every $k$, we have
\[
f^\#(a_k)=|f'(a_k)|=|F'(a_k)e^{g(a_k)}|\ \forall k\in\N
\]
and therefore
\[
f^\#(a_{n_j})\ge |F'(a_{n_j})|\ \forall j\in\N.
\]
Thus it suffices to show that
\[
\lim_{k\to_infty} |F'(a_k)|=+\infty
\]
in order to deduce that $f$ is not Brody.
Now
\[
|F'(a_n)|
=\frac{1}{|a_n|}\left|\prod_{k=1}^{k\ne n}
\left(1-\frac{a_n}{a_k}\right)\right|.
\]
Now
\[
\lim_{n\to\infty}
\frac{1}{|a_n|}
\left|1-\frac{a_n}{a_1}\right|=\frac{1}{|a_1|}.
\]
By the first claim
\[
\left|\prod_{k>n}
\left(1-\frac{a_n}{a_k}\right)\right|\ge C
\]
and by the second claim for each $K>0$ there is a number $N\in\N$
such that
\[
\left|\prod_{1<k<n}
\left(1-\frac{a_n}{a_k}\right)\right|\ge K
\]
for all $n\ge N$.

Combined these assertions show that for each $K>0$ there is a number
$N$
such that
\[
|F'(a_n)|\ge\frac{C}{2|a_0|}K
\]
for all $n\ge N$. Thus $\limsup|f^\#(a_n)|=+\infty$ and $f$ is not Brody.
\end{proof}

\begin{corollary}
Let $\zeta:\R^+\to\R^+$ be an unbounded monotone increasing function.

Then there exists an effective reduced divisor $D$ such that
\begin{enumerate}
\item
For every $r\in\R^+$ the inequality $\deg(D_r)\le\zeta(r)$ holds
where $D_r$ denotes the restriction of $D$ to the open disc with
radius $r$.
\item
There is no Brody entire function $f$ for which $D$ is the zero
divisor.
\end{enumerate}
\end{corollary}

In the language of Nevanlinna theory:
\begin{corollary}
Let $\zeta:\R^+\to\R^+$ be an unbounded monotone increasing function.
Then there exists an effective reduced divisor $D$ such that
\begin{enumerate}
\item
For every $r\in\R^+$ the inequality $N(r,D)\le\zeta(r)$ holds.
\item
There is no Brody entire function $f$ for which $D$ is the zero
divisor.
\end{enumerate}
\end{corollary}
\section{Growth conditions and characteristic function}
Our goal is to show that no bound on the characteristic function
$T_f(r)$ forces an entire function $f$ to be Brody except when this
bound is strong enough to force $f$ to be a polynomial.

In view of the preceding section the crucial point is to verify that
for every such bound there exists an entire function $f$ fulfilling
this condition on $T_f(r)$ and fulfilling simultaneously the
condition of theorem~\ref{div-slow-gro}.

Before stating the result of this section we recall
some basic notions.

For an entire function $f$ with $f(0)\ne 0$ we may define the
characteristic function as $T_f(r)=m_{1/f}(r)+N_{1/f}(r)$
with
\[
m_{1/f}(r)=\int_0^{2\pi}\log^+\frac{1}{|f(re^{i\theta})|}d\theta
\]
and
\[
N_{1/f}(r)=\sum_{a\in\C}\ord_a(f)\log^+\left|\frac{r}{a}\right|.
\]

Now we can state the result:

\begin{theorem}
Let $\rho:[1,\infty[\to]0,\infty[$ be a continuous increasing
function.
\begin{enumerate}
\item
If $\liminf_{t\to\infty}\frac{\rho(t)}{\log t}<\infty$, then every
entire function with $T_f(r)\le\rho(r)$ 
for all $r\ge 1$ must be a polynomial.
\item
If $\liminf_{t\to\infty}\frac{\rho(t)}{\log t}=\infty$, then there
exists an entire function which is not Brody and such that
$T_f(r)\le\rho(r)\ \forall r \ge 1$.
\end{enumerate}
\end{theorem}

\begin{proof}
$(i)$.
If there is a constant $C>0$ such that
\[
\liminf_{t\to\infty}\frac{\rho(t)}{\log t}< C
\]
then
\[
\liminf_{t\to\infty}\frac{N_{f-a}(t)}{\log t}< C
\]
for every entire function $f$ with $T_f(r)\le \rho(r)$
and for all $a\in\C$. It follows that each fiber $f^{-1}\{a\}$
has cardinality $\le C$ and that $f$ is a polynomial of degree
bounded by $C$.

$(2)$.
As a first preparation, we observe that
\[
\prod_k\frac{2^{2k+1}-1}{2^{2k+1}}
=\prod_k\left(1-\frac{1}{2^{2k+1}}\right)
\ge 1-\sum_k \frac{1}{2^{2k+1}}=1-\frac{2}{3}=\frac{1}{3}.
\]

We will construct $f$ as follows: We will choose a sequence $c_k\in\C$
with $|c_{k+1}|>4|c_k|\ge 1$ for all $k$, then define
\[
P_k(z)=3\prod_k\left(\frac{z}{c_k}-1\right).
\]
The conditions on the $c_k$ ensure that $P_k$ converges locally
uniformly to an entire function $f$ whose zeroes are precisely
the points $c_k$.

Let $k\in\N$ and let $z$ be a complex number with
$\frac{1}{2}|c_{k+1}|\le|z|\le 2|c_k|$.
Then
\[
\left|\frac{z}{c_j}-1\right|\ge 1
\]
for all $j\le k$.
For $j>k$ then conditions $|c_{l+1}|\ge 4|c_l|$ imply
\[
\left|\frac{c_j}{z}\right|\ge 2^{2(j-k)+1}
\]
which in turn implies
\[
|P_j(z)|\ge 1
\]
for all $j$ and all $z$ 
with $\frac{1}{2}|c_{k+1}|\le|z|\le 2|c_k|$.
As a consequence,
\[
m_{1/f}(r)=0
\]
and hence
\[
T_f(r)=N_{1/f}(r)
\]
for all $r\in [\frac{1}{2}|c_{k+1}|, 2|c_k|]$ and all $k\in\N$.

Therefore 
\[
T_f(r)=\sum_{j\le k}\log\left|\frac{r}{c_j}\right|\le k\log r
\]
for all such $r$.

Similarily one obtains $T_f(r)=0$ for $r\le\frac{1}{2}|c_1|$.

Let us now fix $k\in\N$ and consider $r\in
[\frac{1}{2}|c_k|, 2|c_k|]$.

Since $T_f(r)$ is increasing, we have
\[
T_f(r)\le T_f(2|c_k|)=N_{1/f}(2|c_k|)=
\sum_{j\le k}\log\left|\frac{2c_k}{c_j}\right|\le k\log(2|c_k|)
\le k\log(4r).
\]

Summarizing, we have shown that $T_f(r)\le k\log(4r)$
for all $r$ with $r\le 2|c_k|$.

Thus it suffices to choose the $c_k$ such that
$k\log(4r)\le \rho(r)$ for all $r\in[2|c_{k-1}|,2|c_k|]$.

This is possible: We assumed 
 $\liminf_{t\to\infty}\frac{\rho(t)}{\log t}=\infty$,
hence for each $k\in\N$ there is a constant $R_k$ such that
$k\log(4r)\le\rho(r)$ for all $r\ge R_k$. Now it suffices to choose
the $c_k$ such that $2|c_{k-1}|\ge R_k$ (in addition to the other
conditions $|c_k|\ge 4|c_{k-1}|\ge 1$).

Thus we have established:
{\em We can choose a sequence $c_k$ in $\C\setminus\{0\}$ such that
$f(z)=3\prod_k\left(\frac{z}{c_k}-1\right)$ converges to an entire
function $f$ with $T_f(r)\le\rho(r)$ for all $r\ge 1$.}

Finally, we note that in our construction we choose
the $c_k$ such that $|c_{k+1}|\ge 4|c_k|$.
Furthermore we may choose the $c_k$ such that
the set
\[
\left\{\frac{c_k}{|c_k|}:k\in\N\right\}
\]
is dense in $S^1=\{z\in\C:|z|=1\}$.
Then theorem~\ref{div-slow-gro}
implies such an entire function $f$ ist not Brody.
\end{proof}

\section{The divisor $D=\{k^2:k\in\N\}$}
\begin{proposition}
Let $a_k=k^2$. Then there is a Brody function with divisor
$\{a_k:k\in\N\}$.
\end{proposition}

\begin{proof}
We pose
\[
f(z)=f_D(z)=\prod_k\left( 1 - \frac{z}{k^2} \right)
\]
It is easily verified that this is convergent.
An explicit calculation shows:
\[
f'(a_k)=f'(k^2)=
\lim_{n\to\infty}\frac{1}{2k^2}\frac{(n-1+k)\cdots n}%
{(n+1)\cdots(n+k)}=\frac{1}{2k^2}.
\]
In particular, $|f'(a)|\le\frac{1}{2}$ for all $a\in D$.

Now let $z\not\in|D|$. We will deduce an estimate
for
\[
\sum_k \frac{1}{|z-a_k|},
\]
since
\[
\frac{f'(z)}{f(z)}=\sum_k \frac{1}{z-a_k}.
\]
First we note: If a subset of $\C$ contains at least $m$ numbers
in $|D|=\{k^2:k\in\N\}$, then its diameter must be at least $m^2-1$.
Let $b_k$ be a renumbering of the points in $|D|$ such that
$|b_{k+1}-z|\ge|b_k-z|$ for all $k$.
If $r=|b_k-z|$ for some $k\in\N$, then 
\[
b_1,\ldots,b_k\in \overline{D_r(z)}=\{w\in\C:|z-w|\le r\}
\]
and therefore $r\ge k^2-1$.
Thus
\[
\frac{1}{|b_k-z|}\le\frac{1}{k^2-1}
\]
for all $k>1$.
We deduce:
{\em For every $z\not\in|D|$ there is a point $b\in|D|$ with
$d(z,D)=|z-b|$ and
\begin{multline}
\left|\left( \sum_k \frac{1}{z-a_k}\right) - \frac{1}{z-b}\right|
=\left| \sum_{k\ge 2} \frac{1}{z-b_k}\right|
\le\sum_{k=2}^\infty \frac{1}{k^2-1}=\\
=\frac{1}{2}\sum_{k=2}^\infty \left(\frac{1}{k-1}-\frac{1}{k+1}\right)
=\frac{5}{12}
\end{multline}
}
Similarily
\[
\sum_{k\in\N\setminus\{m\}}\frac{1}{|a_m-a_k|}
\le\frac{5}{12}
\]
for $m\in\N$.

We are now ready to show that $f^{\#}$ is bounded. 
We start with the case $d(z,D)\ge 1$.
If $d(z,|D|)\ge 1$, then
\begin{align*}
f^{\#}(z)&=\left|\frac{f'(z)}{f(z)}\right| h(f(z)) \\
&\le (1-5/12)h(f(z)) \le  (1-5/12)\frac{1}{2}= \frac{7}{24}
\end{align*}

It remains to discuss the case $d(z,D)\le 1$. 
Let $b\in D$
such that $d(z,D)=d(z,b)=|z-b|$.

As a preparation we discuss 
\[
P(z)=\prod_{a\in D\setminus\{b\}}\frac{a-z}{a-b}
=\prod_{a\in D\setminus\{b\}}\left(1+\frac{b-z}{a-b}\right).
\]
Using $|b-z|\le 1$ we obtain the following bound:
\[
|P(z)|\le
\prod_{a\in D\setminus\{b\}}\left(1+\frac{1}{|a-b|}\right)
\le\exp\left(\sum_{a\in D\setminus\{b\}}\frac{1}{|a-b|}\right)
\le\exp\left(\frac{5}{12}\right).
\]

We know:
\[
f'(b)=-\frac{1}{b}\prod_{a\in D\setminus\{b\}}\left(1-\frac{b}{a}\right)
\]
and
\[
f(z)=\left(1-\frac{z}{b}\right)
\prod_{a\in D\setminus\{b\}}\left(1-\frac{z}{a}\right).
\]
Combined these two equations yield:
\begin{align*}
&(z-b)\frac{f'(b)}{f(z)}
\prod_{a\in D\setminus\{b\}}\left(1-\frac{b}{a}\right)=1\\
\iff &
(z-b)\frac{f'(b)}{f(z)}P(z)=1
\end{align*}

Hence
\begin{align*}
f'(z)&=\frac{f'(z)}{f(z)}f(z)\\
&=\frac{f'(z)}{f(z)}(z-b)f'(b)P(z)\\
&=\left(\frac{1}{b-z}+\sum_{a\in D\setminus\{b\}}\frac{1}{a-z}\right)
(z-b)f'(b)P(z)\\
&=-f'(b)P(z)+\left(\sum_{a\in D\setminus\{b\}}\frac{1}{a-z}\right)(z-b)f'(b)P(z).
\end{align*}
It follows that
\[
|f'(z)|\le \frac{1}{2}\exp\left(\frac{5}{12}\right)+
\frac{5}{12}1\cdot \frac{1}{2}\exp\left(\frac{5}{12}\right)
\]
for all $z\in \C\setminus D$ with $d(z,D)\le 1$.
\end{proof}

\section{Discussion}
Using the special cases of entire functions studied above
we see that the class of entire functions which are Brody is
not closed neither under addition nor under multiplication:
the entire functions $z$, $e^z+1$ and $ze^z$  are all Brody,
but $ze^z+z$ is not, although
\[
ze^z+z= z(e^z+1)=\left(ze^z\right)+z.
\]

We see also that the Brody condition is neither closed nor open
nor complex:

For $(s,t)\in\C^2$ let us consider
\[
f_{s,t}(z)=se^z+\frac{z}{tz-1}
\]
By proposition~\ref{case-1} the function $f_{s,t}$ is Brody iff 
\[
(s,t)\in\{(x,y)\in\C^2:x=0 \text{ or }y\ne 0\}
\]
which is neither a closed nor an open set.

Moreover, proposition~\ref{ex-re} provides an example of a family of entire
functions depending holomorphically on a complex parameter $\lambda$
such that the function is Brody if and only if $\lambda$ is real.

All this properties are in stark contrast to the situation for
Brody curves with values in abelian varieties.
If $A$ is an abelian variety, its universal covering is isomorphic
to some $\C^g$. Since every holomorphic map from $\C$ to $A$ lifts
to a holomorphic map with values in the universal covering of $A$,
the classical theorem of Liouville implies that an entire curve
with values in $A$ is Brody if and only if it can be lifted
to an affine-linear map from $\C$ to $\C^g$.
As a consequence one obtains:
\begin{itemize}
\item
If $f,g:\C\to A$ are Brody curves, so is $z\mapsto f(z)+g(z)$.
\item
If $f_t:\C\to A$ is a family of entire curves depending 
holomorphically on a parameter $t\in P$ where $P$ is a complex
manifold, then the set of all $t\in P$ for which $f_t$ is Brody
forms a closed complex analytic subset of $P$.
\item
An entire curve $f:\C\to A$ is Brody if and only if
\[
\limsup \frac{\log T_f(r)}{\log r} \le 2.
\]
\end{itemize}


\begin{thebibliography}{Bla}

\bibitem{B}
Brody, R.: 
Compact manifolds and hyperbolicity.
\sl T.A.M.S. \bf 235\rm, 213--219 (1978)

\bibitem{CH}
Clunie, J.; Hayman, W.K.:
The spherical derivative of integral and meromrophic functions.
\sl Comm.~Math.~Helv.\rm, 117-148, \bf 40 \rm ( 1966)

\bibitem{L}
Lehto, O.: 
The spherical derivative of meromorphic functions in the neighbourhood
of an isolated singularity.
\sl Comm. Math.~Helv. {\bf 33}\rm, 196--205 (1959)

\bibitem{LV}
Lehto, O.; Virtanen, K.I.:
On the behaviour of meromorphic functions in the neighbourhood of an
isolated
singularity.
\sl Ann. Acad. Sci. Fenn. Ser. A. I.\rm, \bf 240\rm (1957)

\bibitem{T}
Tsukamoto, M.:
A packing problem for holomorphic curves.
arXiv:math.CV/0605353 (2006)

\end{thebibliography}
\end{document}